\documentclass[12pt,letterpaper, leqno]{article}
\usepackage{graphics,graphicx}
\def\pf {{\bf Proof:} \ }
\def\endpf{$\|$ \bigskip}
\def\reals{\hbox{\sl I\kern-.18em R \kern-.3em}}

\newtheorem{theorem}{Theorem} 
\newtheorem{proposition}{Proposition}
\newtheorem{lemma}{Lemma}

\def\pf {{\bf Proof:} \ }
\def\endpf{$\|$ \bigskip}
\def\cK{{\cal K}}

\def\cA{{\cal A}}
\def\cX{{\cal X}}
\def\cU{{\cal U}}
\def\cB{{\cal B}}
\def\cF{{\cal F}}
\def\cD{{\cal D}}

\def\bA{{\bf A}}

\def\bF{{\bf F}}

\def\np{\noindent} 

\begin{document}
\title {The Markov Theorem for transverse knots}
\author {Nancy C. Wrinkle \\ wrinkle@math.uga.edu}
\date{\today}
\maketitle

\begin{abstract} \noindent 
Let $\xi$ be the standard contact structure in oriented $\reals^3 = (\rho, \theta, z)$ given as the kernel of the $1$-form $\alpha = \rho^2 d\theta + dz$.  A transverse knot is a knot that is transverse to the planes of this contact structure. In this paper we prove the Markov Theorem for transverse knots, which states that two transverse closed braids that are isotopic as transverse knots are also isotopic as transverse braids.  The methods of the proof are based on Birman and Menasco's proof of the Markov Theorem in their recent paper \cite{BM02}, modified to the transverse setting.  The modification is straightforward until we get to the special case of preferred longitudes, where we need some new machinery.  We use techniques from earlier work by the author with Birman \cite{BW00},  by Birman and Menasco (\cite{BM4},  for example), and by Cromwell \cite{Cr95}.  
\end{abstract}

\tableofcontents

\section{Introduction}
\label{section:intro}
Let $\xi$ be the standard contact structure in oriented $\reals^3 = (\rho, \theta, z)$ given as the kernel of the $1$-form $\alpha = \rho^2 d\theta + dz$.  A transverse knot is a knot that is transverse to the planes of this contact structure (although knots can be transverse to any contact structure defined on a $3$-manifold in which the knots are embedded, in this paper we are considering only the standard contact structure on $\reals^3$).  In particular, if we parametrize the knot by $(\rho(t), \theta(t), z(t))$, then the knot is transverse if $(\rho(t))^2 > - \frac{z'(t)}{\theta'(t)}$ for all $t$.  Note that we have made a choice of the sign of $\alpha$: we will assume throughout that $\alpha > 0$.   
In this paper we prove the Markov Theorem for transverse knots:\\

\np {\bf Transverse Markov Theorem (TMT): Let $X_1$ and $X_2$ be closed transverse braids in standard contact $\reals^3$, with the same braid axis, and let $X_1$ and $X_2$ be transversely isotopic as transverse knots. In particular, $X_1$ and $X_2$ have the same topological knot type $\cX$ and the same self-linking number.  Then $X_2$ may be transversely obtained from $X_1$ by transverse braid isotopy and a finite number of transverse, that is, positive, stabilizations and destabilizations.} \\

The methods of the proof are based on Birman and Menasco's proof of the Markov Theorem in their recent paper \cite{BM02}, modified to the transverse setting.  We use techniques from \cite{BW00} and \cite{BM02} and develop new methods from Cromwell's paper, \cite{Cr95}.  \\

The essence of the Birman-Menasco proof of the Markov Theorem is showing that the isotopy of $X_1$ and $X_2$ can be divided into two isotopies: the first is taking the connected sum $X_1 \bigoplus U_i$ of $X_1$ with several copies of the unknot, and the second is an isotopy of  $X_1 \bigoplus U_i$ across its Seifert surface to $X_2$.  The modification of the Birman-Menasco proof goes easily for the first isotopy, and the difficulty lies in the second isotopy.  This special case is proved in Section \ref{section:II}, using techniques based on braid foliations followed by an argument based on Cromwell's work with arc-presentations. \\

The paper is organized as follows.  In Section \ref{section:background}, we will include some background on braid foliations, (c.f. \cite{BM02}).  In Section \ref{section:lemmas}, we cover the first isotopy in the TMT.  In Section \ref{section:II} we examine the geometry of a special surface, the annulus bounded by a transverse closed braid and its transverse preferred longitude.  This surface is a subsurface of the Seifert surface for one of the boundary braids.  Then we use our knowledge of the foliation of the annulus to modify the second isotopy of the Markov Theorem to the transverse setting.  In the last section, we put together the results of the previous sections to modify Birman and Menasco's proof of the Markov Theorem for the case of transverse knots.  \\

{\bf Remark:} It is surprising that before this year, no proofs of the Transverse Markov Theorem appeared.  Bennequin proved a transverse version of the Alexander Theorem in \cite{Be}, and other authors (\cite{FT} and others) have defined a series of Reidemeister-type moves on Legendrian and transverse knots.  However, the only proof of an equivalence theorem is in Swiatkowski's paper \cite{Sw}, in which he proves a Reidemeister theorem for non-intersecting plane fronts.  The main result of this paper was proven simultaneously and independently by Orevkov and Shevchisin (http://xxx.lanl.gov/abs/math.GT/0112207).  They used contact geometric methods, and posted their work a few weeks earlier than the appearance of this paper.\\

{\bf Acknowledgements:} This paper is based on the dissertation written in partial satisfaction of the requirements for the PhD at Columbia University.  The author is grateful to Bill Menasco for many helpful conversations.  Thanks to Joan Birman, Bill Menasco, and Elizabeth Finklestein for allowing the reproduction of some of the figures from their papers.   The borrowed figures are cited in their captions.  Special thanks to Joan Birman for the suggestion of this thesis problem and for her many hours of advising and patience  throughout the years in which this work took shape.\\

\section{Background on braid foliations.}
\label{section:background}
Although we are proving a result about closed braids, all of our work is guided by the surface bounded by the closed braids, and a certain foliation of this surface.  Birman and Menasco developed a detailed analysis of this foliation for use in their series of papers with the common title "Studying links via closed braids" (for example, see \cite{BM4}).  There is a review of their methods in \cite{BF}, which includes proofs of many of the results described below.  In this section we will review the terminology and methods from these papers, following the background given in \cite{BM02}. \\

Let $K$ be a closed braid with braid axis $\bA$, representing the knot type $\cK$.  Let $\cal{H}$ be the open book decomposition of $\reals^3$ with open half-planes $H_{\theta}$ through constant $\theta$ as leaves and the axis $\bA$ as its binding.  As a closed braid, $K$ intersects each $H_{\theta}$ transversely and is oriented in the direction of increasing $\theta$.  We assume that $K$ bounds a surface $\cF$ of maximum Euler characteristic.  The foliation that we will study in this paper is the intersection of the surface $\cF$ with the leaves of $\cal{H}$.  This foliation is not the characteristic foliation of the surface that is given by the integration of the line field of the surface intersected with the planes of the contact structure.  \\

  After the surface has been modified (the modification is described in \cite{BF}), it admits a special foliation studied by Birman and Menasco, which we will use in this paper.  We may assume that the surface has the following properties.  The intersections of the axis $\bA$ with the surface $\cF$ are transverse and finite in number.  Each intersection of $\bA$ with $\cF$ is called a {\em vertex} and has a disk neighborhood which is radially foliated.  In a neighborhood of $K$ the foliation of the surface is radial.  All but finitely many of the fibers $H_{\theta}$ meet $\cF$ transversely and those which do not are tangent to the surface at exactly one point of the surface and the fiber.  Each point of tangency is a hyperbolic singularity.  Such a surface is called a {\em Markov surface}.  \\

\begin{figure}[htpb]
\centerline{\includegraphics[height=2.75in]{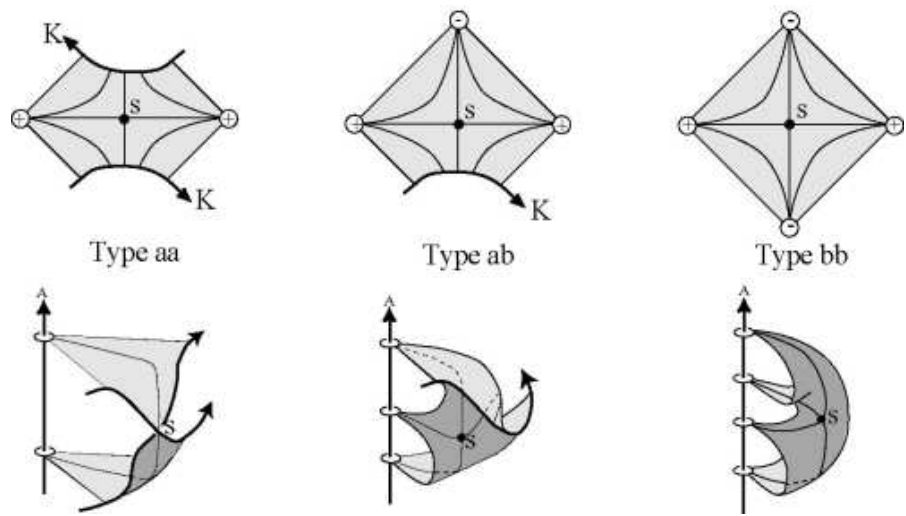}}
\caption {The three tile types and their embeddings. (\cite{BM02}, \cite{BF})}
\label{figure:tiles}
\end{figure}

A {\em singular leaf} in the foliation is one which contains a point of tangency with some $H_{\theta}$.  All other leaves are {\em nonsingular}.  Each singular leaf contains exactly one saddle singularity.  Nonsingular leaves do not have both endpoints on $K$ and are divided into two  types: $a$-arcs, which have one endpoint on $K$ and one on $\bA$, and $b$-arcs, which have both ends on $\bA$.  The foliation may be used to decompose the surface into a union of foliated $2$-cells, called {\em tiles}, each containing one singularity.  The vertices of the tiles are the vertices of the foliation, and the edges are $b$-arcs or $a$-arcs.  The tiles fall into three types, $aa$, $ab$, and $bb$, according to the type of nonsingular leaves coming together before the singularity.  See Figure \ref{figure:tiles}, which shows the tiles and their embeddings.  \\

The vertices and singularities of the foliation have signs determined by the orientation of the surface.  If the axis pierces the surface from the negative (resp. positive), side then the vertex is positive (resp. negative).  If the outward pointing normal of the surface at a singularity points in the direction of increasing (resp. decreasing) $\theta$, then the singularity is positive (resp. negative).  The {\em valence} of a vertex $v$ is the number of singular leaves with an endpoint at $v$.\\

Next we describe isotopies of $K$ and of the interior of the surface bounded by $K$ that are guided by the tiles of the foliation of the surface $\cF$.  These moves all simplify the foliation either by reducing the number of vertices and singularities or by reducing the valence of the vertices.  \\

\begin{figure}[htpb]
\centerline{\includegraphics[height=2.5in]{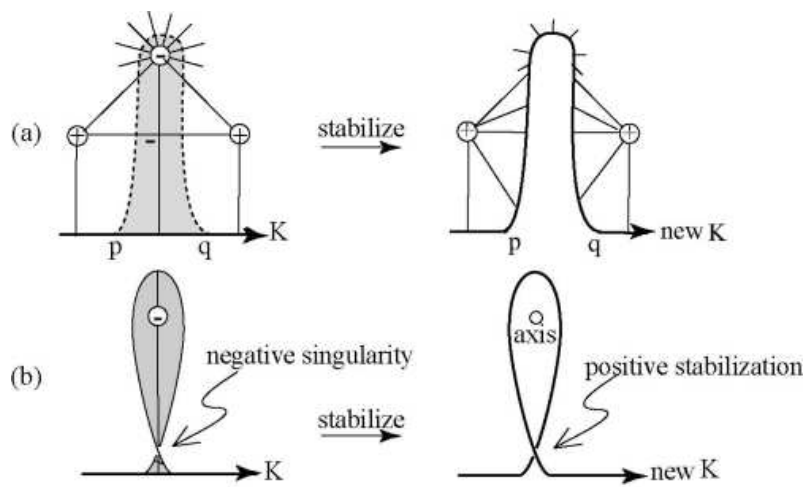}}
\caption {Stabilizing along an $ab$-tile. \cite{BM02}}
\label{figure:abstab}
\end{figure}

\np {\bf Stabilization along $ab$-tiles:} Figure \ref{figure:abstab} shows stabilization of $K$ along an $ab$-tile.  It consists of pushing $K$ across a disk  neighborhood of a singular leaf in a tile to get rid of the singularity and a negative vertex $v$.  Looking at the effect on $K$ of stabilization along an $ab$-tile, we see that it is the familiar Markov move on a closed braid with n strands that makes it into a closed braid with n+1 strands.   Therefore stabilization increases the braid index of $K$ by 1.  The effects on the foliation of the stabilization along an $ab$-tile are to delete the $ab$-tile and the negative vertex $v$, change all $ab$-tiles that have $v$ as a vertex to $aa$-tiles and change all $bb$-tiles that have $v$ as a vertex to $ab$-tiles.  The sign of the added trivial loop is the opposite of the sign of the singularity in the tile (Figure \ref{figure:abstab}(b)). \\

\begin{figure}[htpb]
\centerline{\includegraphics[height=1.5in]{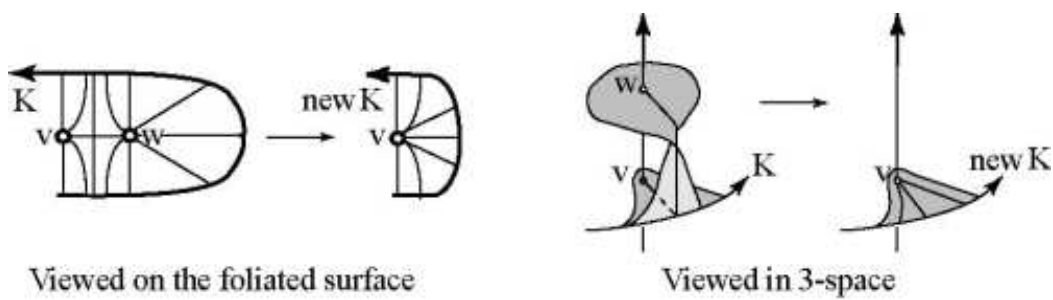}}
\caption {Destabilizing along an $aa$-tile. \cite{BM02}}
\label{figure:aastab}
\end{figure}

\np {\bf Destabilization along $aa$-tiles:} If the foliation has a vertex of valence $1$,  then two edges of a single tile are identified in the foliation.  The tile is an $aa$-tile  (Figure \ref{figure:tiles}) with the pieces of the boundary braid identified in a trivial loop, and a neighborhood of the valence 1 vertex is a radially foliated disc.  See Figure \ref{figure:aastab}.  The $aa$-destabilization along an $aa$-tile removes the trivial loop, using the Markov move that takes a closed (n+1)-braid to a closed n-braid.  The effects on the foliation of destabilization along an $aa$-tile are to delete the $aa$-tile and the valence 1 vertex.  It does not change the foliation in the complement of the $aa$-tile.  The sign of the trivial loop is the same as the sign of the singularity in the $aa$-tile.    \\

Note that stabilization along an $ab$-tile and destabilization along an $aa$-tile both reduce the number of vertices and singularities in the foliation.  Therefore, although stabilization and destabilization are inverse moves on $K$, they are not inverse moves on the foliation of the surface. \\

\begin{figure}[htpb]
\centerline{\includegraphics[height=2.5in]{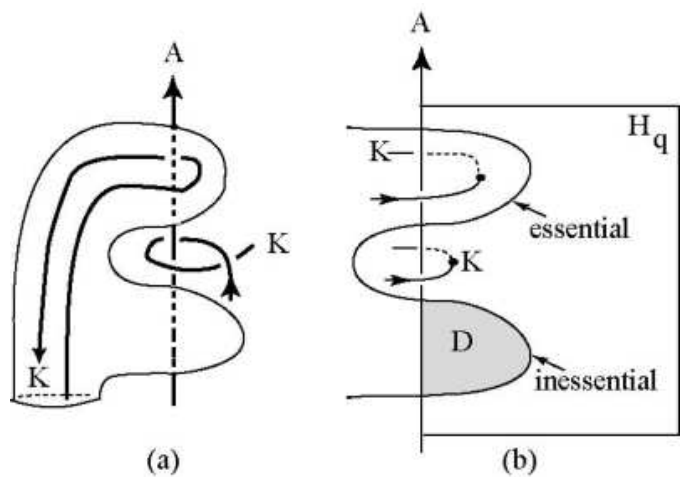}}
\caption {Essential and inessential $b$-arcs. \cite{BM02}}
\label{figure:inessential}
\end{figure}

Another isotopy of the interior of the surface, called {\em removing an inessential $b$-arc}, allows us to simplify the foliation further.  It reduces the number of vertices and singularities in the foliation each by 2.  Each $b$-arc in a fiber $H_{\theta}$ separates $H_{\theta}$ into two disks.  If $K$ pierces both of these disks, then the $b$-arc is called {\em essential}.  Otherwise the $b$-arc is {\em inessential}.  See Figure \ref{figure:inessential}.  An alternative way of looking at inessential $b$-arcs is to note that a $b$-arc is inessential if its endpoints are adjacent along the braid axis.\\

\begin{figure}[htpb]
\centerline{\includegraphics[height=2.5in]{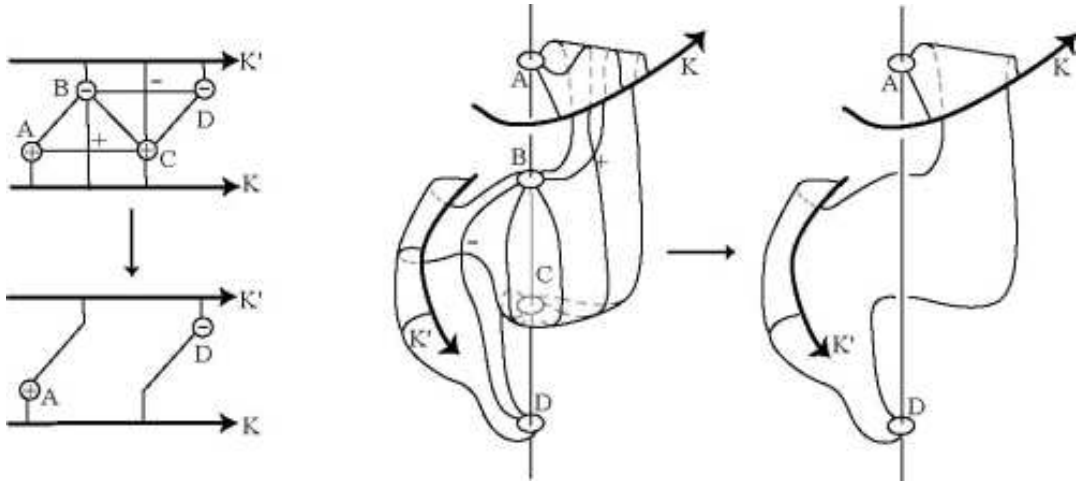}}
\caption {Removing an inessential $b$-arc from an embedded annulus.}
\label{figure:inessential square}
\end{figure}

Another view of a piece of an embedded surface with an inessential $b$-arc is given in Figure \ref{figure:inessential square}, which will be an important example later in this paper.  The surface is an annulus that is bounded by two closed braids $K$ and $K'$, and the inessential $b$-arc is shown coming out of the page between vertices $B$ and $C$.  The modified surface with the inessential $b$-arc and surrounding singularities and vertices removed is shown on the right.\\ 

\np \begin{lemma} All $b$-arcs may be assumed to be essential. \end{lemma}
\np \pf Any inessential $b$-arc may be removed via an isotopy of the surface that removes two vertices and two singularities from the foliation.  See Figure \ref{figure:inessential square}. For details, see \cite{BF}.  \endpf

\begin{figure}[htpb]
\centerline{\includegraphics[height=1in]{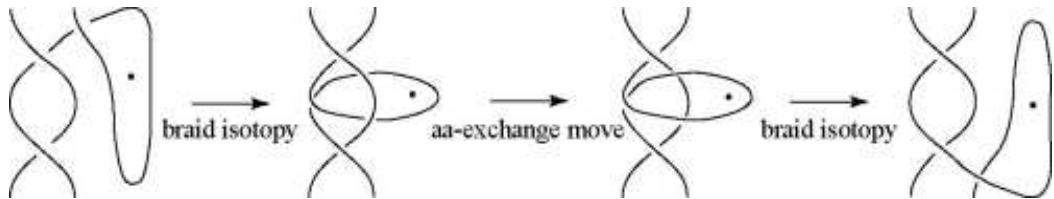}}
\caption {An $aa$-exchange move on $K$. }
\label{figure:aaexch}
\end{figure}

\begin{figure}[htpb]
\centerline{\includegraphics[height=1.5in]{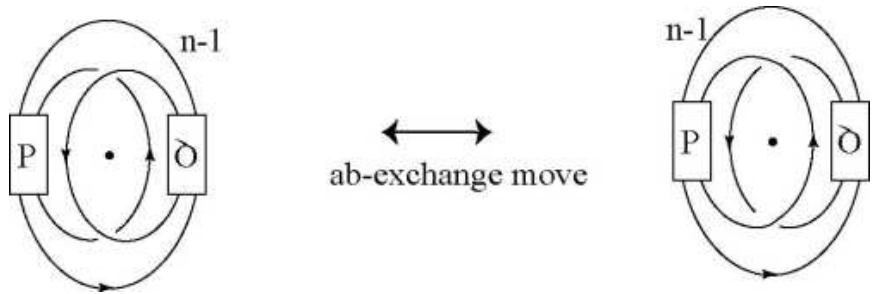}}
\caption {Performing an $ab$-exchange move on a closed $n$-braid.}
\label{figure:exchange}
\end{figure}

\begin{figure}[htpb]
\centerline{\includegraphics[height=1.5in]{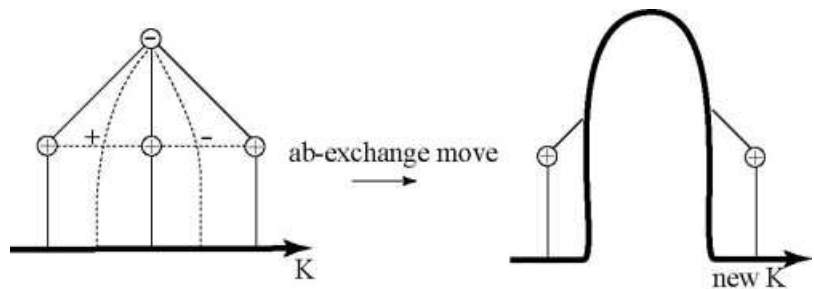}}
\caption {Effects of an $ab$-exchange move on the foliation of the surface.}
\label{figure:abtiles}
\end{figure}

{\bf Exchange moves:}  An {\em exchange move} is a modification of the braid and the surface that is admitted when a vertex has valence 2 and the signs of the singularities in the two adjacent tiles is $(\pm, \mp)$.  There are three types of exchange moves, classified according to the tiles of the vertex of valence 2: $aa$, $ab$, and $bb$.  The effects of the $aa$- and $ab$-exchange moves on the boundary braid are shown in Figures \ref{figure:aaexch} and \ref{figure:exchange}.  After the $aa$- or $ab$-exchange moves, the foliation of the modified surface has two fewer vertices and singularities than the original foliation.  For example, see Figure \ref{figure:abtiles}.\\ 

\begin{figure}[htpb]
\centerline{\includegraphics[height=1.5in]{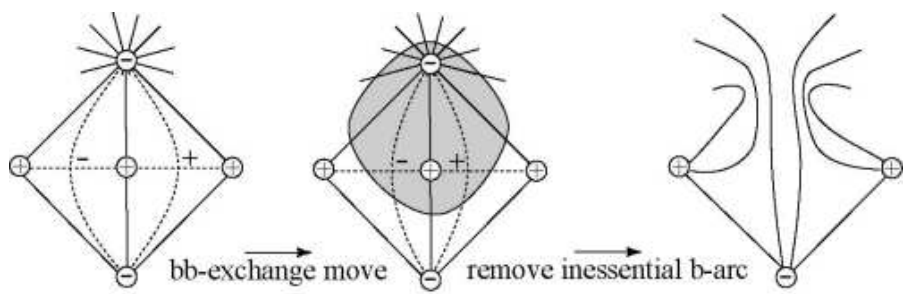}}
\caption {A $bb$-exchange move.}
\label{figure:bbexch}
\end{figure}

The $bb$-exchange move has the same effect on the boundary braid as the $ab$-exchange move.    After the $bb$-exchange move, the number of vertices and singularities of the foliation are unchanged, but there is an inessential $b$-arc in the foliation.  See Figure \ref{figure:inessential}(a), which shows part of a braid and its surface.  If we move the upper loop below the lower loop, the $b$-arc that was essential is now inessential.  This appearance of an inessential $b$-arc via a $bb$-exchange move is marked in Figure \ref{figure:bbexch} by a gray disc.  After removing the inessential $b$-arc, the number of vertices and singularities in the foliation are both reduced by 2.  \\

Although the exchange move changes the foliation, it does not change the braid index, as is clear from Figures \ref{figure:exchange} and \ref{figure:aaexch}.  \\

\begin{figure}[htpb]
\centerline{\includegraphics[height=1.75in]{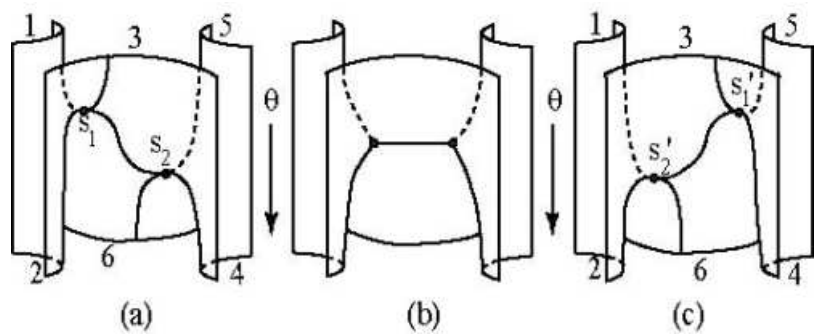}}
\caption {A change in foliation on an interior disc of $\cF$.\cite{BF}}
\label{figure:foliation}
\end{figure}

\begin{figure}[htpb]
\centerline{\includegraphics[height=1.5in]{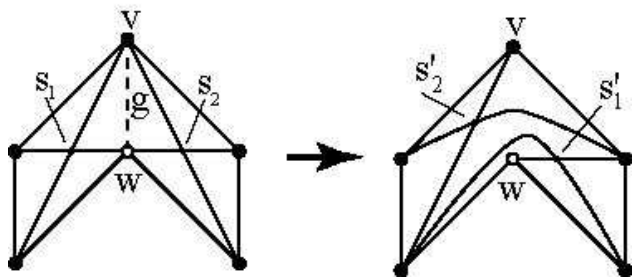}}
\caption {A change in foliation reduces the valence of vertices $v$ and $w$.\cite{BF}}
\label{figure:foliated tiles}
\end{figure}

Another modification of the foliated Markov surface $\cF$ that we will use is called simply a {\em change in foliation}.  In a change in foliation, adjacent, same signed singularities are moved past each other in the $\theta$ direction.  See Figure \ref{figure:foliation}.  This change takes place entirely on the interior of the surface, without changing the boundary braid at all.  Since it does not change the boundary braid, it will not affect the transverse knot type of the boundary braid.  The immediate effect of the change in foliation is to reduce the valence of two of the vertices in the foliation.  See, for example, the change in foliation that takes place along the two $bb$-tiles in Figure \ref{figure:foliated tiles}.  \\

In the following work, along with combinations of the moves described above, we will be using braid isotopies in the complement of the braid axis.  We call such an isotopy simply a {\em braid isotopy}. 

\section{Part One of the TMT}
\label{section:lemmas}

As mentioned in the Introduction, the proof of the Markov Theorem consists of dividing the given braid isotopy into two pieces.  The first piece of the isotopy is taking the connected sum of a braid $K$ with $n$ copies of the unknot.  In this section we give the transverse modification of this isotopy.  In order to preserve the Bennequin number of the braid $K$, when we take connected sums with copies of the unknot, we assume that each of the unknots has Bennequin number $-1$.  \\

Throughout the paper, we will be using stabilization, destabilization, exchange moves, and braid isotopies.  The following lemma, then, is essential. \\

\begin{lemma} 
\label{lemma:BW}
{\bf (see page 335 of \cite{BW00})}  If a transverse closed braid is modified by one of the following isotopies, then the isotopy can be replaced by a transverse isotopy:
\begin{enumerate}
\item Braid isotopy.
\item Positive stabilization or destabilization.
\item An exchange move. 
\end{enumerate}
\end{lemma}

\np \begin{lemma}
\label{lemma:unknot}
{\bf (The unknot lemma)}:
Let $U$ be a transverse closed braid with knot type the unknot
and $\beta(U) = -1$.  Then $U$ can be changed by a finite number of transverse moves into the standard 1-braid
representative of the unknot.  \end{lemma}

\np \pf The unknot knot type $\cU$ is transversally simple. 
Therefore, by \cite{BW00}, there exists a sequence

$$U \rightarrow U_1 \rightarrow U_2 \rightarrow \cdots U_{n-1} \rightarrow  
U_n$$

\np where $U_n$ is the standard 1-braid representative, and each $U_i\in \cU$ is given by a destabilization, exchange move, or braid isotopy
of the preceding $U_{i-1}$.  \\

Going through the sequence, each time we come to a negative destabilization, we can use an $aa$-exchange move to slide the negative trival loop out of the way, pushing the negative destabilization along the sequence.  (See Figure \ref{figure:aaexch} in Section \ref{section:background}).  In this way, we move all of the $U_i$ that result from negative destabilization to the end of the sequence, and we get a new sequence:

$$U \longrightarrow \{U'_i\}_{i=1}^r \longrightarrow 
\{U^{\star}_j\}_{j=0}^{s-1} \longrightarrow U_n = U_s^{\star}$$

\np where each of the $U'_i$ results from a positive destabilization, exchange
move, or braid isotopy
of the preceding $U'_{i-1}$, and each of the $U^{\star}_j$ result from negative destabilization or braid isotopy of $U^{\star}_{j-1}$. \\

By construction,  $\beta(U'_r)
= \beta(U) = -1$, and in fact by Lemma \ref{lemma:BW} they are transversally isotopic. 
Also by construction,
$\beta(U^{\star}_s) = -1 + 2s$, since every negative
destabilization increases $\beta$ by 2.  However, $\beta(U_n) =
-1$ and $U_n= U^{\star}_s$.  It must be that $-1 = -1 +
2s$ so  $s = 0$ and $U'_r = U_n$.  \endpf

\np \begin{lemma}
\label{lemma:braidconnectsum}
{\bf (The connected sum lemma)}:
Let $B$ be a closed transverse braid with knot type $\cB$ and
let $U$ be a transverse closed braid with knot type the unknot
and $\beta(U) = -1$.  If the braid connected sum of $B$ and
$U$, $B \bigoplus_b U$ has Bennequin number $\beta(B \bigoplus_b U) = \beta(B)$, then $B \bigoplus_b U$ and $B$ are transversely isotopic as closed braids.  
\end{lemma}

\np \pf Holding $B$ fixed, we may transversely modify $U$ as in Lemma \ref{lemma:unknot} until it is the standard 1-braid representative of the unknot.  Call the modified braid $U'$.  Let $\cF$ be the Markov surface for $B$ and $\cD$ be the standard disc bounded by the standard 1-braid representative for the unknot.  Then $B \bigoplus_b U'$ bounds $\cF \bigoplus_b \cD$, which is just $\cF$.  To check that the boundary braids are transversely isotopic we note that by construction, the addition of the standard disc to the surface $\cF$ corresponds to an increase in the braid index of the original braid $B$ but no change in the Bennequin number.  Using the braid index formula for Bennequin number ($\beta(K) = e(K) - n(K)$, where $e(K)$ is the algebraic crossing number of $K$ and $n(K)$ is its braid index,) we conclude that the trivial loop bounding the disc has a positive crossing, so $B$ and $B \bigoplus_b U$ are transversely isotopic as closed braids via a positive destabilization.  \endpf

\np \begin{lemma}
\label{lemma:nbraidconnectsum}
{\bf (The n connected sum lemma)}:
Let $B$ be a braid as in lemma \ref{lemma:braidconnectsum} and
let $\{U_1, U_2, \ldots, U_s\}$ be a finite sequence of
transverse unknots with $\beta (U_i) = -1$  for all $i$.   Let
$B_i$ be defined recursively by $B_0 = B$ and $B_j = B_{j-1}
\bigoplus U_j$ (Here $\bigoplus$ is not the braid connected sum
$\bigoplus_b$ but the connected sum of $B_{j-1}$ and $U_j$ as
knots).  Then
$B_s$ is transversely isotopic to $B$ as closed braids.  
\end{lemma} 

\np \pf The composite braid theorem \cite{BM4} states that given
a composite link type $\cK$, there is a closed $n$-braid
representative $K$ and a sequence of closed $n$-braids
$\{K_i\}$ taking $K$ to some $K'$ given by exchange
moves and braid isotopy,
such that $K'$ is a composite $n$-braid.  Since both
exchange moves and braid isotopy in the complement of the
braid axis are transverse isotopies, this theorem, combined
with Bennequin's transverse Alexander Theorem \cite{Be}, which is stated below for the benefit of the reader, gives
a transverse isotopy between a transverse composite link (as
defined, for example, by the recursive formula above), and a transverse
composite closed braid. \\

\np {\bf Transverse Alexander Theorem (see Theorem 8 of \cite{Be}):} {\em If $K$ is a knot that is transverse to the standard contact structure in $\reals^3$, then it can be transversely isotoped to a transverse closed braid.}\\

Suppose $s=1$.  Then $B_s =B \bigoplus U_1$.  By the transverse
composite braid theorem, $B \bigoplus U_1$ is transversely isotopic
to $B' \bigoplus_b U'_1$, where $B' $ and $
U'_1$ are closed transverse braids and $B'$ is
transversely isotopic to $B$ and $ U'_1$ is
transversely isotopic to $U_1$.  Then by Lemma
\ref{lemma:unknot},  $B' \bigoplus_b U'_1$ is
transversely isotopic to $B'$, which is transversely
isotopic to $B$.  \\

Suppose $s>1$, and consider $B_s = B_{s-1} \bigoplus U_s$.  Again we
transversely make $B_s$ into a transverse composite closed
braid $B_{s-1}' \bigoplus_b U_s'$.  By Lemma
\ref{lemma:unknot}, we know that $B_{s-1}' \bigoplus_b
U_s'$ is transversely isotopic to $B_{s-1}'$. 
Continuing on, we have that $B_{s-1}'$ is transversely
isotopic to $B_{s-2}^{\prime\prime} \bigoplus U_{s-1}'$
where $B_{s-2}^{\prime\prime}$ is transversely isotopic to
$B_{s-2}$, and again we use Lemma \ref{lemma:unknot} to get
from $B_{s-2}^{\prime\prime} \bigoplus U_{s-1}'$ to 
$B_{s-2}^{\prime\prime}$.  Each of these steps may produce a
different transverse closed braid but each $B_{s-j}^{(j)}$ is
transversely equivalent to $B_{s-j}$.  So at the end, 
$B_s$ is transversely equivalent to $B_0^{(s)}$, and $B_0^{(s)}$
is transversely equivalent to $B$.  \endpf

\section{Part Two of the TMT.}
\label{section:II}
\subsection{Geometry of the annulus bounded by transverse preferred longitudes}
\label{subsection:geometry}

Let $K$ be a transverse closed braid and let $K'$ be a preferred longitude for $K$, with $\beta(K) = \beta(K')$.  We may assume that $K'$ is also a transverse closed braid by the transverse Alexander's Theorem.  In this section we focus on the second part of the isotopy of the transverse Markov Theorem, in which we push $K$ across an embedded annulus to $K'$.  In order to proceed we need to understand better the geometry of the annulus bounded by $K$ and $K'$.  \\

Because $K'$ is a preferred longitude for $K$, they bound an embedded annulus $\cal{A}$.  This annulus is a subsurface of a Markov surface for $K$.  In this section, we will extend some of Birman's and Menasco's results about the foliated surface to the transverse setting and use these results to simplify the foliation of $\cal{A}$.   The simplification consists of the transverse moves on the boundary braids $K$ and $K'$ as described in Lemma \ref{lemma:BW}, and changes in foliation on the interior of the surface as described in Section \ref{section:background}.  Our goal for now is to simplify the foliated annulus enough to describe its embedding completely.  Following the methods described in Section \ref{section:background},  we divide our work among the three types of tiles.  \\

{\bf Remark:}  In the annulus the surface will be to the left of $K$, and to the right of $K'$, while all of the moves guided by the surface have been described in Section \ref{section:background} for the case when the surface is to the left of the braid.  The braid being to the right of the surface simply changes the sign of the braid moves being guided by the surface.  For example, even though only positive stabilizations are transverse and the $ab$-tiles that guide these moves are negative for $K$ (as described in Section \ref{section:background}), for $K'$ a positive stabilization is guided by a positive $ab$-tile. \\

We will now describe the foliation that results from these simplifications.  In a later section we will describe the embedding of the annulus, which will allow us to further simplify the foliation.  

\begin{figure}[htpb]
\centerline{\includegraphics[height=3in]{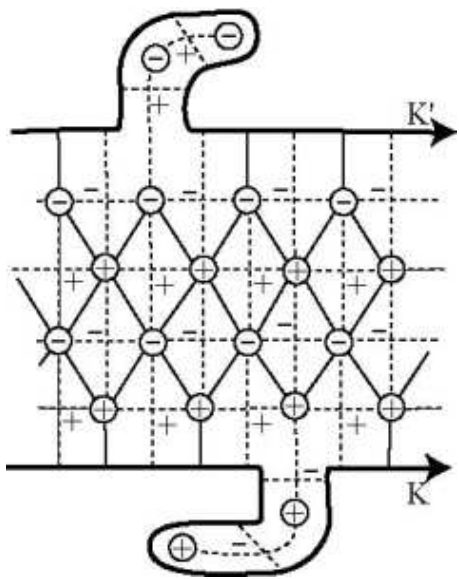}}
\caption {The foliation of the annulus bounded by a transverse closed braid and its transverse preferred longitude.}
\label{figure:full annulus}
\end{figure}

\begin{proposition} 
After some number of transverse stabilizations and destabilizations, braid isotopy, and exchange moves of the transverse closed braid $K$ and its preferred longitude $K'$, and some number of changes of foliation in the embedded annulus bounded by $K$ and $K'$, we may assume that the annulus has a checkerboard foliation as pictured in Figure \ref{figure:full annulus}, along with a tab of negative $aa$-tiles along $K$ and a tab of positive $aa$-tiles along $K'$, each tab bounded by a collection of trivial negative loops of $K$ and $K'$.  \end{proposition}

\pf
We start our with the $aa$-tiles.  Each closed $aa$-tile is bounded by a trivial loop in either $K$ or $K'$, with the sign of the trivial loop of $K$ corresponding to the sign of the tile and the sign of the trivial loop of $K'$ the opposite of the tile's sign. Depending on the sign of the tile, we may be able to transversely (positively) destabilize the braid along the tile to get rid of the vertex and singularity.  We can destabilize $K$ if the tile bounded by $K$ is positive, and $K'$ if the tile bounded by $K'$ is negative.  After we have transversely destabilized $K$ and $K'$ along as many of the $aa$-tiles as possible, we are left with some number of closed negative $aa$-tiles along $K$ and some number of closed positive $aa$-tiles along $K'$.  Using the $aa$-exchange moves described in Section \ref{section:background}, we can push the $aa$-tiles along each boundary until they are all in the same spot, forming a tab from $K$ and $K'$.   Recall Figure \ref{figure:aaexch}.  It will be useful later on to consider the braids $K$ and $K'$ minus these tabs of $aa$-tiles. Although they will not be the same transverse knot type as $K$ and $K'$, they will be the same topological knot type.  \\ 

Along each braid there will be some number of positive and negative $ab$-tiles.  We first eliminate pairs of oppositely signed tiles that share a common $a$-arc and a common $b$-arc by performing $ab$-exchange moves along $K$ or $K'$, as described in Section \ref{section:background}.  Recall Figure \ref{figure:abtiles}.  If there are any $ab$-tiles along $K$ or $K'$ that share a common $b$-arc and have the same sign, then we may reduce the valence of  their common vertices by using a change in foliation.  Recall Figure \ref{figure:foliated tiles}.  Depending on the sign of the tiles, we may stabilize $K$ or $K'$ further along the modified $ab$-tiles.   If there are any other negative $ab$-tiles along $K$, we may transversely remove them by positively stabilizing $K$ along those tiles.  In the same way, we may get rid of any other positive $ab$-tiles along $K'$ by positively stabilizing $K'$ along them, noting that the sign is correct because $K'$ is to the right of the surface.  Recall Figure \ref{figure:abstab}.  We conclude that after some number of $ab$-exchange moves and $aa$-exchange moves, positive stabilizations, and changes in foliation, along $K$ we are left with only positive $ab$-tiles that share an $a$-arc and the tab of negative $aa$-tiles constructed above, and along $K'$ we have negative $ab$-tiles sharing an $a$-arc and the tab of positive $aa$-tiles. \\

To complete our description of the foliated annulus, we need some lemmas and terminology from \cite {BF}:  \\

\np {\bf Definitions:} An {\em interior} vertex is a vertex in the foliation that is not the endpoint of an 
$a$-arc.  {\em Star($v$)} is the subset of the foliation consisting of the vertex $v$, all $b$-arcs with an endpoint $v$, all singular leaves with an endpoint $v$, and the vertices and singularities at the end of these $b$-arcs and singular leaves.  Star($v$) gives us purely combinatorial information about the vertex $v$.  The {\em star order} of a vertex $v$ is the number of edges in Star($v$).   See Figure \ref{figure:star alone}. \\

\begin{lemma} 
\label{lemma:BF3.1}
{\bf (See Lemma 3.1 on page 310 of \cite{BF})}
 Let $v$ be an interior vertex.  Then $star(v)$ contains both positive and negative singularities.
\end{lemma}

\begin{figure}[htpb]
\centerline{\includegraphics[height=2.5in]{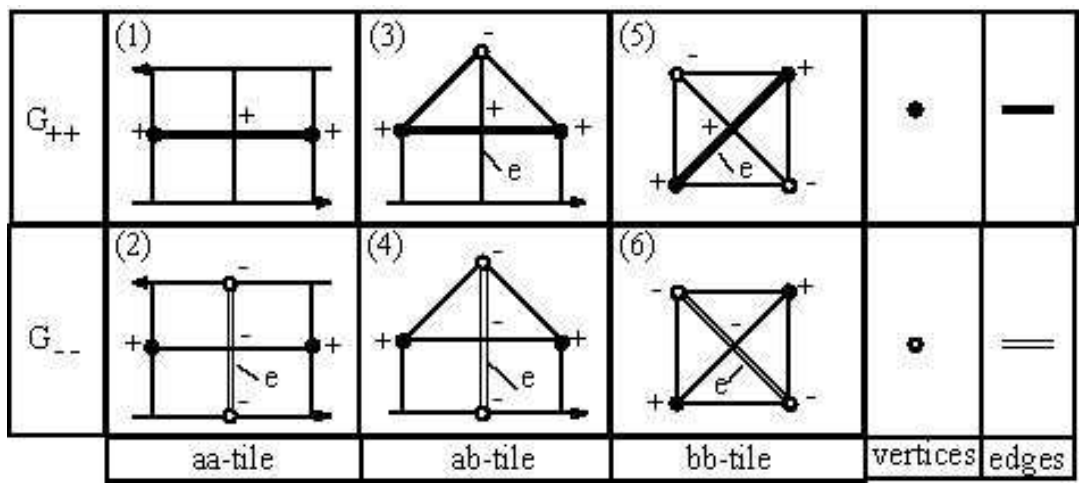}}
\caption {The graphs $G_{++}$ and $G_{--}$ of a foliated surface. \cite{BF}}
\label{figure:graphs}
\end{figure}

\begin{figure}[htpb]
\centerline{\includegraphics[height=2in]{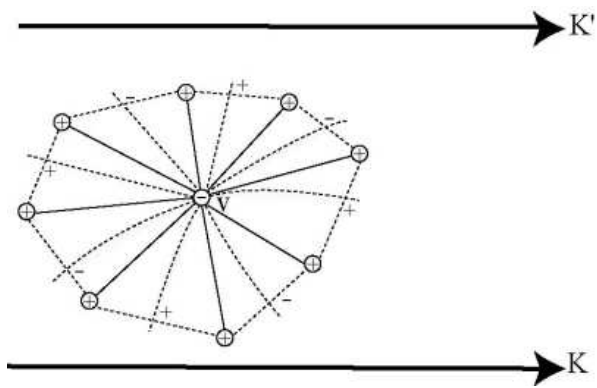}}
\caption {Star($v$), where $v$ is an interior vertex with star order $= 16$.}
\label{figure:star alone}
\end{figure}

\np {\bf Definition:} The {\em graph} $G_{\epsilon, \delta}$ of an embedded foliated surface is a graph whose vertices are all of the vertices of sign $\epsilon$ and whose edges are all of the singular leaves which go through singularities of sign $\delta$.  See Figure \ref{figure:graphs}, which shows some of the possibilities for the graphs $G_{++}$ and $G_{--}$, the only graphs we will need for our purposes.  The figure does not include the graphs through tiles adjacent to the braid $K'$ to the right of the surface.  We can easily obtain those tiles by reflecting the given tiles across the braid, reversing the signs of the vertices and singularities, and swapping the edges and vertices of $G_{++}$ and $G_{--}$.  These graphs first appeared in Bennequin's work (see \cite{Be}).  For example, in the annulus (Figure \ref{figure:full annulus}), the graphs $G_{++}$ and $G_{--}$ form concentric, alternating closed curves in the surface.  \\

The following results hold for embedded surfaces bounded by a braid, but by Lemma \ref{lemma:BW}, they are also true for an embedded surface with transverse boundary.  \\

Let $\vec x = x, \dots, x$, where $x$ is a sign $+$ or $-$. \\

\begin{lemma} 
\label{lemma:BF3.2} {\bf (See Lemma 3.2 on  page 313 of  \cite{BF})}
After some number of changes in foliation, exchange moves, and braid isotopies, we may assume that no interior vertex has star with the sign of its singularities $(\vec +, \vec -)$. \end{lemma}

\begin{lemma} 
\label{lemma:BF3.7} {\bf(See Lemma 3.7 on Page 318 of  \cite{BF})} After some number of changes in foliation, exchange moves, and braid isotopies, we may assume that no interior vertex can be an endpoint of the graph $G_{++}$ or $G_{--}$.\end{lemma}

\begin{lemma} 
\label{lemma:BF3.8i}{\bf (See Lemma 3.8(i) on Page 318 of  \cite{BF})}
After some number of changes in foliation, exchange moves, and braid isotopies, we may assume that  no closed loop of any graph $G_{++}$, $G_{--}$, $G_{+-}$, or $G_{-+}$ of the foliation will bound a disk. \end{lemma}

We complete our description of the foliation of the interior of the surface with the following lemma. \\

\begin{lemma}After some number of changes in foliation, exchange moves, and braid isotopies, we may assume that all interior vertices have valence $2$ in the graphs $G_{++}$ or $G_{--}$.  \end{lemma}

\pf 
We will prove the claim for the case when the vertex is negative.  If a negative interior vertex $v$ has valence 1 in the graph $G_{--}$, then it is an endpoint of the graph.  By Lemma \ref{lemma:BF3.7}, after performing some number of exchange moves, changes in foliation, and braid isotopies, we may assume that $v$ is not the endpoint of $G_{--}$, so does not have valence 1 in the graph. \\

\begin{figure}[htpb]
\centerline{\includegraphics[height=3in, width=6in]{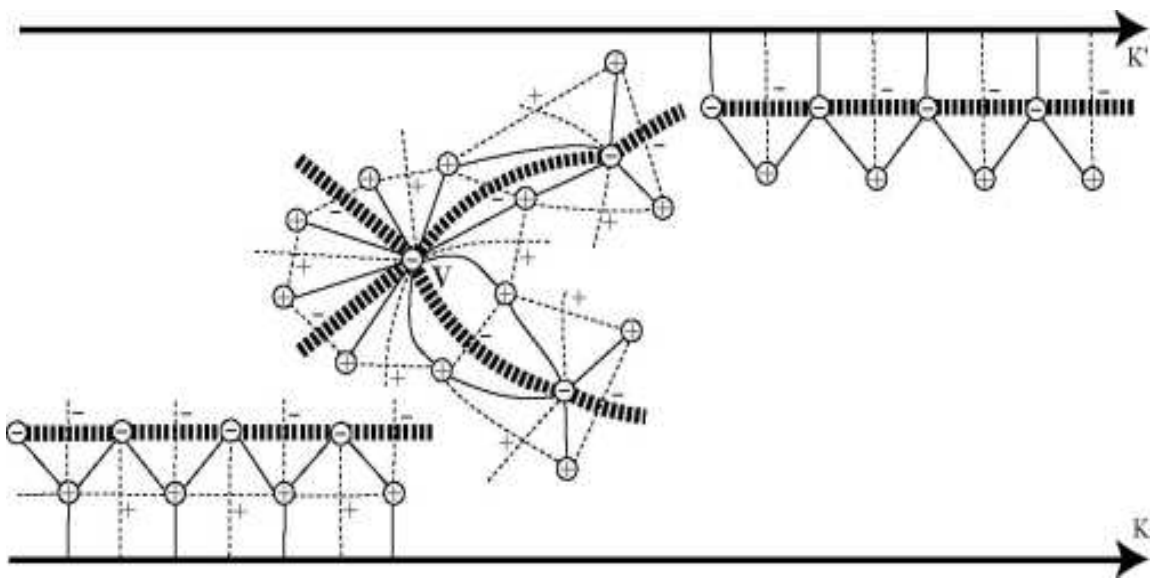}}
\caption {The graph $G_{--}$ of part of the foliated annulus, if the annulus contains an interior vertex $v$ whose valence is greater than 2 in the graph.}
\label{figure:wrong graph}
\end{figure}

If an interior vertex has valence greater than 2 in the graph , then there is some branching of the graph at $v$.  See Figure \ref{figure:wrong graph}.  After some number of exchange moves, changes in foliation, and braid isotopies, the graph $G_{--}$ does not have endpoints on $K$ or $K'$, or at any interior vertex, so it must be a closed graph.  If there is branching, then only one of the branches closes up to form a curve that is parallel to the boundary and all other branches close up as the boundary of discs.  By Lemma \ref{lemma:BF3.8i}, after some number of exchange moves, changes in foliation, and isotopies in the complement fo the axis, we may assume that no closed loop of $G_{--}$ bounds a disc.  Therefore, there is no branching, so no negative interior vertex in the graph has valence greater than 2.    \\

The proof of the lemma is complete when we observe that the graph $G_{++}$, after some number of exchange moves, changes in foliation, and braid isotopies, also has no endpoints on $K$ or $K'$ so the results above about $G_{--}$ and negative interior vertices is also true for $G_{++}$ and positive interior vertices.  \endpf

To conclude that the surface is foliated as in Figure \ref{figure:full annulus}, we perform more changes in foliation, exchange moves, and isotopies in the complement of the axis until we can assume that for each interior vertex $v$, the singularities in $star(v)$ alternate sign (Lemma \ref{lemma:BF3.2}).  Since by the above lemma we know that each interior vertex has exactly two positive or two negative vertices, this is enough to conclude that the valence of each interior vertex in the surface is 4, and our foliation has the claimed form.  The proof of the proposition is complete.  \endpf

 An example of a foliated annulus with all interior vertices having  valence $2$ in the graphs is pictured in Figure \ref{figure:full annulus}.  Although we didn't know at the beginning of this section how many positive $ab$-tiles and negative $ab$-tiles there are in the foliated annulus, from the configuration of the $bb$-tiles we see that there must be the same number of each.  \\

We may use our results about the foliation to conclude that the braid index of $K$ must be the same as the braid index of $K'$, modulo the {\em tab} of $aa$-tile bounding trivial loops coming off of each (see Figure \ref{figure:full annulus}).  For in terms of the surface bounded by the braid, the braid index of the boundary braid is the algebraic intersection number of the axis with the surface.  This gives a formula for the braid index of the braid from its foliated surface: $b(K) = $ \#(positive vertices)$ - $ \#(negative vertices).  We already have the Markov surface $\cF$ for $K$, and we may isotope the Markov surface for $K'$ to $\cF-\cA$, the subsurface of $\cF$ that is the complement of the checkerboard annulus.  Then the difference between $b(K)$ and $b(K')$ is the difference \#(positive vertices)- \#(negative vertices) in the annulus $\cA$, which is 0, modulo the vertices in the $aa$-tiles.  We shall see later that in fact there is the same number of positive and negative $aa$-tiles.  \\

\begin{figure}[htpb]
\centerline{\includegraphics[height=2in]{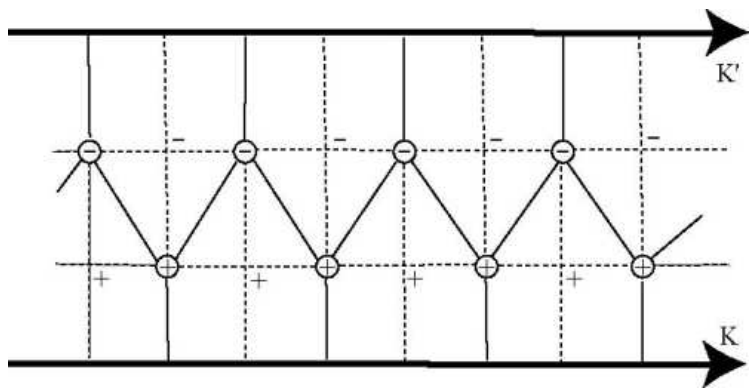}}
\caption {The simplified annulus without trivial loops.}
\label{figure:simple annulus}
\end{figure}

We have shown that, after some number of changes in foliation, exchange moves, transverse stabilizations and destabilizations of $K$ and $K'$, we may assume that the foliation of the annulus they bound has the checkerboard pattern pictured in Figure \ref{figure:full annulus}.  In fact, it suffices to prove the TMT for preferred longitudes for the case when the foliation has an even simpler pattern: a checkerboard with no $bb$-tiles.  See Figure \ref{figure:simple annulus}.  For if we can push $K$ across this simplified annulus to $K'$, then we may inductively work across the rest of an annulus with $bb$-tiles in the pattern of Figure \ref{figure:full annulus}.  \\

It is here that the modification of the Birman-Menasco proof reaches a difficulty.  In the non-transverse setting, we may destabilize $K$ and $K'$ along the $ab$-tiles and stabilize them along the $aa$-tiles, giving us a nonsingularly foliated annulus, then we may push $K$ across the remaining surface to $K'$.  Since none of these moves is transverse, we must attack the foliation differently, by showing that if $K$ and $K'$ are transverse preferred longitudes, then the moves we have already performed leave us with a foliation whose $b$-arcs are all inessential. Recall that there is an embedding of the annulus of Figure \ref{figure:simple annulus} with all of the $b$-arcs inessential: a piece of it is pictured in Figure \ref{figure:inessential square}.  The inessential b-arcs in this figure are facing out of the page, between vertices $B$ and $C$.  When we consider an extension of this embedding through the rest of the annulus, with all of the $b$-arcs inessential, we will be able to remove the vertex $A$ and the vertex above it, and the vertex $D$ and the vertex below it, and continue in this way through the rest of the foliation.  This leaves an annulus which does not intersect the axis at all.\\

\subsection{The arc-presentation of an oriented knot type.}
\label{subsection:arcpres definitions}

Until now, we have considered the annulus as a single embedded surface with two boundary components, $K$ and $K'$.  Now we would like to shift our point of view to the graphs $G_{++}$ and $G_{--}$ of the foliation, and describe the annulus as a union of  subannuli with these graphs as core circles.  In the following, we will use Cromwell's work on arc-presentations to first decompose the annulus into a union of three subannuli.  Then we will prove that this decomposition can always be replaced by a simpler decomposition into a union of two subannuli whose foliation has all inessential $b$-arcs. \\

Arc-presentations of oriented knot types have been studied by Cromwell, whose paper \cite{Cr95}  we closely follow here.\\

Consider the open-book decomposition of the 3-sphere which has open disks $H_{\theta}$ for pages and an unknotted circle (which we can take as the $z$-axis $\cup$ $ \infty$ ) for the binding.  A knot can be embedded in finitely many half-planes $H_{\theta}$ so that it meets each half-plane in a single simple arc.  Such an embedding is called an {\it arc-presentation} of the knot. Cromwell proves that every oriented link type has an arc-presentation.  We will be concerned with arc-presentations of oriented knots. \\  

\begin{figure}[htpb]
\begin{center}
\includegraphics[height=2in]{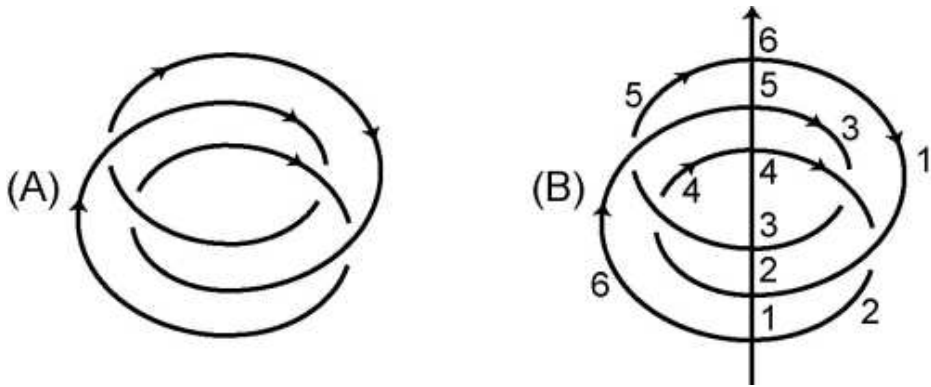}
\end{center}
\caption{(A) The figure 8 knot and (B) its arc-presentation.}
 \label{figure:figure8}
 \end{figure}

As an example, consider the Figure $8$ knot pictured in Figure \ref{figure:figure8}.  The knot diagram pictured in  Figure \ref{figure:figure8} is easily turned into an arc-presentation. There are many ways to obtain an arc-presentation from a given braid diagram; they are detailed in \cite{Cr95} and \cite{Cr96}.\\

An arc-presentation for an oriented knot type $\cal{K}$ is determined (up to equivalence by a sequence of moves described in Section \ref{subsection:cr's construction}) by two sets of data: the cyclic list of points of intersection of the knot with the $z$-axis, ordered with respect to increasing $z$ and listed as one traverses the knot once, and the cyclic list of half planes each arc appears in, ordered with respect to $\theta$ and again listed as one traverses the knot.  In the example in Figure \ref{figure:figure8}, the arc-presentation is given by the permutations $(1,5,3,6,2,4)$ and $(6,3,5,1,4,2)$ ordering the intersections and half-planes, respectively.   The data describing an arc-presentation in this way corresponds directly to the  ordered lists of vertices and singularities obtained from traversing one component of the graph $G_{++}$ or $G_{--}$ of the embedded simplified annulus that we described in Section \ref{subsection:geometry}.  In fact, the data that describe $G_{++}$ and $G_{--}$ determine $G_{++}$ and $G_{--}$ as an arc-presentation of $\cal{K}$.  \\  

\np \begin{lemma}
\label{lemma:graph to arcpres} After an isotopy of the interior of the annulus that fixes the boundary braids, we may assume that the graph $G_{++}$ (and $G_{--}$) of the foliated annulus $\cal{A}$ is an arc-presentation of the oriented knot type $\cal{K}$. \end{lemma}

\np \pf   We prove the claim for the graph $G_{++}$.  The proof for $G_{--}$ is identical.  There are finitely many edges in the graph $G_{++}$.  Because each edge is a singular leaf in the foliation of the surface by the half-planes $H_{\theta}$, each edge $e_{ij}$ is isotopic to a single simple arc which is embedded in a single half-plane.  If $\{v_i,  v_j\}$ are the endpoints of the edge (i.e., the vertices of the graph), then the singularity $s_k$ between $v_i,$ and $ v_j$ gives the $\theta$-coordinate of the half-plane $H_{ij} = H_{\theta_{s_k}}$ in which $e_{ij}$ may be embedded. Each singularity appears in a separate half-plane.  Therefore $G_{++}$ is an arc-presentation of some knot type. Let $ \hat{K}$ be the braid $K$ minus the trivial negative loops bounding the tab of $aa$-tiles.  The braid $\hat{K}$ has the same topological knot type $\cal{K}$ as $K$, though a different transverse knot type.  To show that the knot type of the arc-presentation is that of $\hat{K}$, note that the foliation of the surface between $G_{++}$ and $\hat{K}$ is nonsingular.  Therefore there is an isotopy $h: S^1 \times [0,1] \rightarrow S^1$ taking $\hat{K}$ to $G_{++}$ which is a braid isotopy on $[0,1)$.  \endpf

\begin{figure}[h]
\begin{center}
 \includegraphics[height=3in]{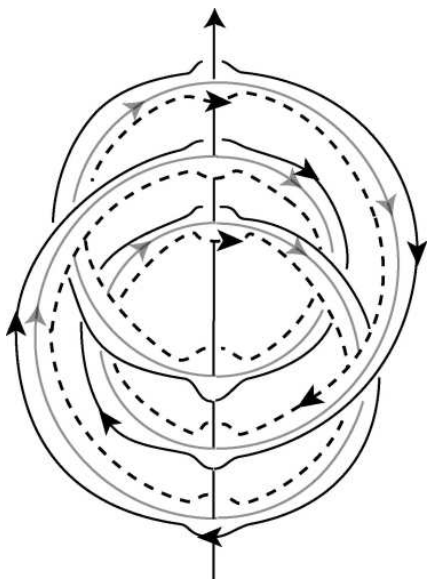}
\end{center}
 \caption {An embedding of an arc-presentation (gray) of the Figure 8 knot, its associated braid (dashed), and its associated anti-braid (solid).}
 \label{figure:embedded arc-pres}
 \end{figure} 

By Lemma \ref{lemma:graph to arcpres}, for the rest of this section we may assume that our surface has been isotoped so that the graphs $G_{++}$ and $G_{--}$ {\em} are arc-presentations for $\cal{K}$.  We can associate a braid to an arc-presentation in a simple way.  At each point of intersection of the arc-presentation with the axis, push the arc-presentation off of the axis in the direction that is normal to the plane spanned by the tangent along the (oriented) arc and the axis.  There is a choice of two directions; choose the normal that makes the positive frame with the tangent and the axis.  Then orientation of the pushed off curve will be positive with respect to $\theta$.  Traveling along the rest of the curve we may push slightly off of the arc-presentation in the same direction, so that the resulting curve is a braid.  This pushoff would be the inverse of the isotopy $h$ used in the proof of Lemma \ref{lemma:graph to arcpres}.  If we push off in the other direction from the axis and along the arc-presentation, we have an anti-braid.  An {\em anti-braid} is a knot that satisfies all of the definitions of a closed braid except that it is increasing with respect to $-\theta$ rather than $\theta$.  We call these knots the braid and anti-braid {\em associated} to the arc-presentation.  See Figure \ref{figure:embedded arc-pres}. (There is an "anti-Markov" theorem for anti-braids, which uses anti-braid isotopy and the Markov moves of type II and III along with stabilization and destabilization, all of which are oriented in the opposite direction from the moves of the Markov Theorem for braids).  \\

\begin{figure}[htpb]
\centerline{\includegraphics[height=2in]{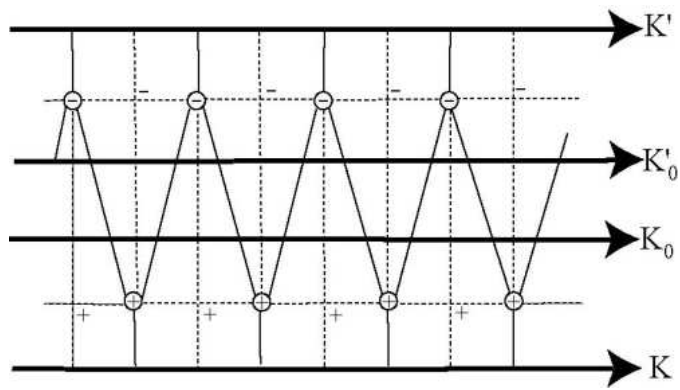}}
\caption {The decomposition of the foliated annulus into three subannuli.}
\label{figure:decomp}
\end{figure}

We use the associated braids and anti-braids of the graphs $G_{++}$ and $G_{--}$ to decompose the annulus into a union of three subannuli.   Given the arc-presentations $G_{++}$ and $G_{--}$ lying in the annulus $\cal{A}$, we can push them off of the axis towards the boundary of the annulus to form two closed braids.  By construction, we may take these braids to be $K$ and $K'$ minus their tabs of negative $aa$-tiles and positive $aa$-tiles, respectively.  We can also push the arc-presentations off of the axis in the other direction into the annulus to obtain two closed anti-braids, $K_0$ and $K'_0$. In this way, we have a decomposition of $\cal{A}$ into three subannuli: one with core circle $G_{++}$, whose boundary components are the braid $K$ and an anti-braid $K_0$, one with core circle $G_{--}$, whose boundary components are the braid  $K'$ and an anti-braid $K'_0$, and a nonsingularly foliated center annulus, whose boundary components are the anti-braids $K_0$ and $K'_0$.  These annuli intersect along the boundary anti-braids. See Figure \ref{figure:decomp}. \\

Now we wish to consider the isotopy of the anti-braid $K_0$ to the anti-braid $K'_0$ that is guided by the nonsingular subannulus of $\cal{A}$ between $K_0$ and $K'_0$.  Recall Figure \ref{figure:decomp}.   The facts that the annulus is nonsingular between $-K_0$ and $-K'_0$ (so $-K_0$ and $-K'_0$ have the same braid index,) implies that the isotopy between $-K_0$ and $-K_0$ is a braid isotopy in the complement of the braid axis.  The isotopy of $-K_0$ corresponds to an isotopy of $K_0$, so we conclude that the anti-braids $K_0$ and $K'_0$ are also related by anti-braid isotopy in the complement of the axis.\\

\subsection{Arc-equivalence and braid isotopy.}
\label{subsection:cr's construction}

In \cite{Cr95} Cromwell defined a set of moves that relate two different arc-presentations of the same oriented knot type. See the thick black curves in Figures \ref{figure:i}, \ref{figure:ii}, and \ref{figure:iii}.  He proved that any two arc-presentations of the same oriented knot type are related by a finite sequence of the moves from this set.  \\

\np {\bf Cromwell's Markov Theorem for Arc-Presentations: (\cite{Cr95}, page 45)} {\em Any two arc-presentations of (oriented link type) $\cal{L}$ are related by a finite sequence of the moves described in Figures \ref{figure:i}, \ref{figure:ii}, and \ref{figure:iii}.}  \\

There is a similar move of type $i(b)$ (Figure \ref{figure:i},) with oppositely signed crossing.  The sign of the crossing depends on the $\theta$-order of the planes involved. Note that the type $ii$ move (Figure \ref{figure:ii},) is possible only when the arcs coming into and out of the two vertices are not in alternating $\theta$-planes.  When looking down the axis at the arc-presentation, this condition is the same as requiring that the arcs through these vertices, though touching at the axis, do not cross.  \\
 
Because $K_0$ and $K'_0$ are braid isotopic, we may choose an isotopy of $K_0$ across the annulus to $K'_0$ that leaves $K$ and $K'$ fixed.   Because $G_{++}$ is isotopic to $G_{--}$ across a nonsingular annulus,  we may also choose an arc-equivalence of $G_{++}$ and $G_{--}$ that fixes $K$ and $K'$. \\

\np \begin{proposition} \label{proposition:iii} To preserve $K$ and $K'$, the arc-equivalence taking $G_{++}$ to $G_{--}$ through the annulus $\cA$ consists only of Cromwell moves of type $iii$.  \end{proposition}

\np \begin{pf}  Notice that the type $iii$ move shown in Figure \ref{figure:iii} is a shift along adjacent $\theta$-planes of the arc-presentation that just changes the order of the singularities without changing the foliation itself.  \\

\begin{figure}[h]
\begin{center}
\includegraphics[height=2.5in]{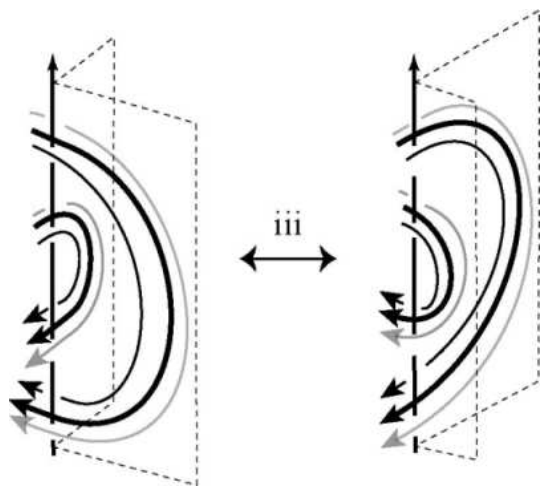}
\end{center}
\caption{A type $iii$ Cromwell move and its effects on the braid and anti-braid.}
 \label{figure:iii}
 \end{figure}

The essence of this proof is to consider the arc-equivalence under the restrictions that the annulus $\cA$ must remain embedded, and that the boundary components $K$ and $K'$ of $\cA$ must be preserved.  Therefore, in Figures \ref{figure:i}, \ref{figure:ii}, and \ref{figure:iii} below, we are looking at changes in the arc-presentation, braid, and anti-braid and considering the embedding of the subannulus bounded by the braid and anti-braid.  Although the braids associated to $G_{++}$ and $G_{--}$ are actually $K$ and $K'$ minus their negative trivial loops, our examination concerns moves that we may assume are supported away from the loops.  Therefore, fixing the braids associated to $G_{++}$ and $G_{--}$ will fix the boundary braids $K$ and $K'$.  We will look at each of the three possible Cromwell moves of the arc-equivalence from this point of view.\\ 

A move of type $i(a)$ or $i(b)$ on $G_{++}$ will stabilize or destabilize either $K$ or $K_0$.   See Figure \ref{figure:i}.  If the braid isotopy of $K_0$ induced type $i$ moves on $G_{++}$ then there would be an induced stabilization or destabilization of $K$.   However, we assumed that the isotopy of $G_{++}$ to $G_{--}$ fixed $K$.  We conclude that the arc-equivalence does not include moves of type $i(a)$ or $i(b)$.\\

\begin{figure}[h]
 \centerline{\includegraphics[height=1.65in]{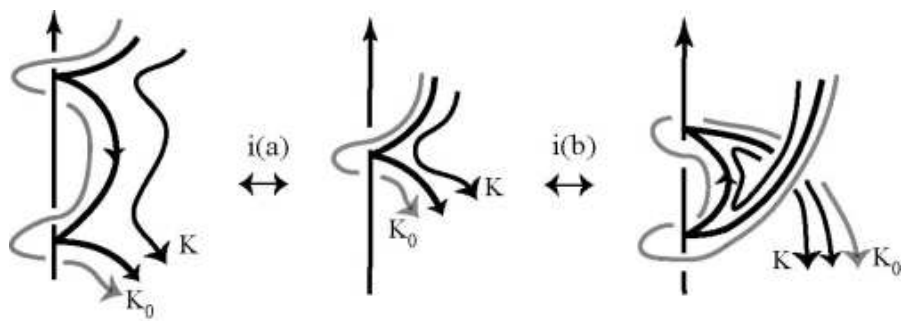}}
 \caption {Cromwell type $i$ moves and the associated isotopy on the braid and anti-braid.}
 \label{figure:i}
 \end{figure}

In Figure \ref{figure:ii}, we see that a move of type $ii$ on two vertices of the arc-presentation corresponds to an exchange move of the braid and a vertical isotopy of a piece of the anti-braid, or to an exchange move of the anti-braid and a vertical isotopy of a piece of the braid, or to a vertical isotopy of pieces of both the braid and the anti-braid, depending on the $\theta$-order and orientation of the involved arcs in the arc-presentation.    Figure \ref{figure:ii}(A) shows an exchange move of the anti-braid (gray) and a vertical isotopy of two pieces of the braid (black).  Figure \ref{figure:ii}(B) shows a vertical isotopy of two pieces of the anti-braid (gray), and a vertical isotopy of two pieces of the braid (black).  \\

Neither $K$ nor $K_0$ admit exchange moves, so if the arc-equivalence includes a type $ii$ move then there is a corresponding vertical isotopy of pieces of $K$ and $K_0$, as in Figure \ref{figure:ii}(B).  However, such a vertical isotopy would cause the subannulus bounded by $K$ and $K_0$ to intersect itself.  We conclude that the arc-equivalence does not include a type $ii$ move.   \\

\begin{figure}[h]
\begin{center}
\includegraphics [height=1.5in]{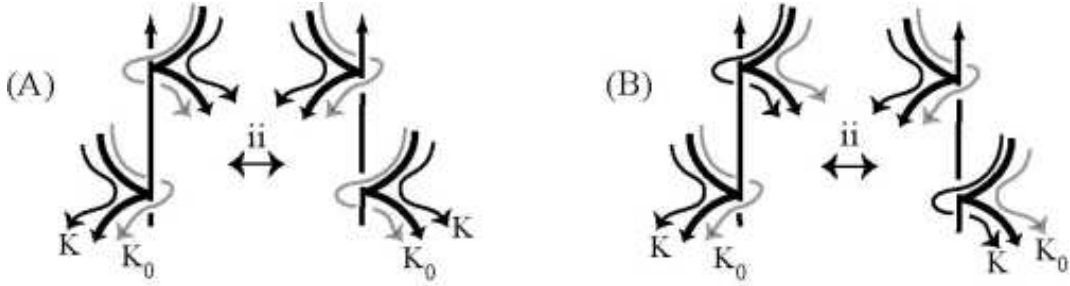}
\end{center}
\caption{A Cromwell type $ii$ move and its effects on the braid and anti-braid.}
 \label{figure:ii}
 \end{figure} 

Ruling out moves of type $i$ and $ii$, we are left only with the moves of type $iii$.  \end{pf} 

\np \begin{proposition} The anti-braid isotopy takes inessential $b$-arcs to inessential $b$-arcs.
\end{proposition}

\np \begin{pf}  By Proposition  \ref{proposition:iii}, we know that only type $iii$ moves occur in the arc-equivalence of $G_{++}$ to $G_{--}$, in order to preserve the boundary and the embedding of the annulus. As a move of type $iii$ on the arc-presentation corresponds to a local $\theta$ shift of the involved strands of the braid and anti-braid,  the isotopy of $K_0$ to $K'_0$ consists simply of shifts of the anti-braid strands through adjacent planes (that is, adjacent in the arc-presentation).  Such an isotopy fixes the order of the vertices in the foliation, so any inessential $b$-arcs are taken to inessential $b$-arcs.   
\end{pf}

Since the isotopy of $K_0$ to $K'_0$ does not change the anti-braid word, we can identify them as a single anti-braid $K_0$.    We get a simpler decomposition of the annulus from this identification of the two anti-braids, just like the one in Figure \ref{figure:decomp}, except that the single anti-braid $K_0$ is the common boundary of the two subannuli.  In this identification, we see that the $b$-arcs connecting the vertices of the arc-presentations traced the trivial shift of $K_0$ to $K'_0$, so are inessential. Therefore, the embedding of the annulus bounded by $K$ and $K'$ is like that pictured in Figure \ref{figure:inessential square}, with all  inessential $b$-arcs.  \\

Using this embedding, we may now prove the following proposition:\\

\np \begin{proposition} {\bf (TMT for preferred longitudes)}  Let $K$ and $K'$ be closed transverse braids in standard contact $\reals^3$ and let $K'$ be a preferred longitude for $K$.  Then $K'$ may be transversely obtained from $K$ by transverse braid isotopy and a finite number of transverse stabilizations and destabilizations. 
\end{proposition}

\pf Let us modify the surface to remove all of the inessential $b$-arcs, as in Figure \ref{figure:inessential square}.  This modification removes all of the singularities and vertices from the foliation, except for the tab of negative $aa$-tiles bounded by negative trivial loops of $K$ and the tab of positive $aa$-tiles bounded by negative trivial loops of $K'$.  We can now push $K$ transversely across the nonsingular annulus and positively stabilize it across the positive $aa$-tiles bounded by $K'$, and we can push $K'$ transversely across the annulus and positively stabilize it across the negative $aa$-tiles bounded by $K$.  This leaves us with a nonsingular foliation of $\cal{A}$. We conclude that there is a transverse braid isotopy that takes $K$ across the nonsingular annulus to $K'$.    \endpf

\section{Proof of the transverse Markov Theorem.}
\label{section:main theorem}

The following proof follows that of Birman and Menasco in \cite{BM02}, using the material of the previous sections to adapt it to the case when the braids are transverse.\\

\np \begin{theorem} \bf{(The transverse Markov Theorem)}:  Let $X_1$ and $X_2$ be closed transverse braids in standard contact $\reals^3$, with the same braid axis, and let $X_1$ and $X_2$ be transversely isotopic as transverse knots. In particular, $X_1$ and $X_2$ have the same topological knot type $\cX$ and the same self-linking number.  Then $X_2$ may be transversely obtained from $X_1$ by transverse braid isotopy and a finite number of transverse stabilizations and destabilizations. 
\end{theorem}

\np \begin{pf}  We will use the notation $X_1 \equiv X_2$ to mean that $X_1$ and $X_2$ are transversely isotopic as closed braids.\\

By hypothesis, we are given closed transverse braids $X_1$ and $X_2$ which represent the same oriented transverse knot type in standard contact $\reals^3$.  We must prove that $X_1 \equiv X_2$.  We may assume without loss of generality that $X_1$ and $X_2$ are situated in distinct half-spaces (so that they are geometrically unlinked), with $X_2$ far above $X_1$.  To prove that $X_1 \equiv X_2$, we will construct a series of knots $X'_2, X^{\prime \prime}_2, X'_3, X_3$, where $X_2$ is a preferred longitude for $X_3$ and $X_3 = X_1 \bigoplus U_i$ is the transverse braid-connected sum of $X_1$ with some number of copies of the unknot described in Section \ref{section:lemmas}.  By the TMT for preferred longitudes, it will follow that $X_1 \equiv X_3$.  By the $n$ connected sum Lemma, it will follow that $X_3 \equiv X_2$. \\

To begin the construction choose a Markov surface $\bF_2$ for $X_2$ (See Section \ref{section:background}) and a transverse preferred longitude $X'_2 \subset \bF_2$ for $X_2$.  We will assume that $X'_2$ lies in a collar neighborhood of $X_2$, chosen to be small enough that $X'_2$ is also a transverse closed braid.  We will also assume that $X'_2$ lies below $X_2$ everywhere except for little hooks where it is forced to travel over a strand of $X_2$, because in general $X_2 \cup X'_2$ is not a split link.  Applying the Unlinking Lemma (below), we construct a knot $X^{\prime \prime}_2$ which is transversally isotopic to $X'_2$ and geometrically unlinked from $X_2$.  The transverse isotopy is a push of $X'_2$ across $k$ disjoint discs $D_1, \dots, D_k$, replacing each $\alpha_i \subset \partial D_i$ by $\beta_i = \partial D_i \ \alpha_i$.  By modifying the subarcs $\beta_i$ of $X^{\prime \prime}_2$ a little bit, if necessary, we may assume that the $\beta_i$ are transverse to every fiber $H_{\theta}$ of the braid foliation and are transverse to the contact structure.  Therefore, $X_2, X'_2$, and $X^{\prime \prime}_2$ are all closed transverse braids.  \\

\np \begin{lemma}
\label{lemma:unlinking}
{\bf (The unlinking lemma)}:
Let $X$ and $X'$ be transverse closed braids that are
not separated by any $2$-plane in $\reals^3$.  Then there is a
transverse isotopy of $X'$ to a new transverse closed braid
$X^{\prime\prime}$ via pushing $X'$ across a set of disjoint
embedded transverse disks, and $X^{\prime\prime}$ and $X$ are
separated by a plane.
\end{lemma}

\np \pf  We can consider $X' \coprod X$ as a transverse
closed braid with more than one component, if we position each
braid so that they are in general position with respect
to the projection $\pi$ onto a  plane $P$ which is
perpendicular to the $z$-axis. (By general position we
mean that the crossings of the different components are
finite in number and produce only finitely many topologically
transverse double points).  A vertical shift  of $X'$ in the
complement of the
$z$-axis can be made transverse, and we
call the result of this isotopy $X''$.  Following the vertical
transverse isotopy of $X'$, the braids $X$ and $X''$ are
linked at finitely many double points $p_1, p_2, \ldots, p_m
\in X$ and
$p'_1, p'_2, \ldots, p'_m \in X''$.  \\

\begin{figure}[h]
\begin{center}
\includegraphics[height=1.5in]{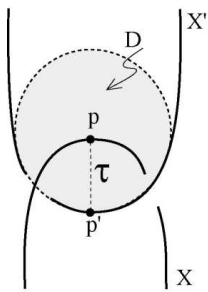}
\end{center}
\caption{A vertical isotopy in the complement of $X$ unlinks the knots.}
 \label{figure:unlinking}
 \end{figure}

See Figure \ref{figure:unlinking}. Let $\tau_i$ be a vertical arc from $p'_i$ to $p_i$.  Let $D_i$
be the disk whose radius is $\tau_i$, part of whose
boundary is a neighborhood of $p_i' \subset X''$, and which is
pierced once by $X$ at the point $p_i$.  The isotopy of $X''$
that will unlink it from $X$ consists of pushing $X''
\cap \partial D_i$ across $D_i$ to $\partial D_i - (X''
\cap \partial D_i)$.  In the complement of $X$ this is
simply a braid isotopy so
can be made transverse, even if the disks $D_i$ are not
transverse.  Repeating the isotopy for $i = 1, 2, \ldots, m$,
we get a transverse closed braid which is geometrically
unlinked from
$X$ and $X'$ and which is transversely equivalent to $X'$.  \endpf

It will be convenient to assume that $X_2$ (resp. $X^{\prime \prime}_2$ and $X_1$) lie in the half-space $\reals^3_+$ (resp. $\reals^3_-$), where the half-spaces are separated by a plane that we refer to as $\reals^2$, also that $X'_2$ lies in $\reals^3_-$ everywhere except for the $k$ hooks where it passes up and over $X_2$, intersecting $\reals^2$ twice.  We think then of each disc $D_i, i = 1, \dots, k$ as a tall thin semi-circular disc which is divided by $\reals^2$ into a semi-circular disc $D_{i,+} \subset \reals^3_+$ and a rectangular disc $D_{i, -} \subset \reals^3_-$, so that $D_{i,-} \cap X^{\prime \prime}_2 = \beta_i$, and $\beta_i = D_i \cap X^{\prime \prime}_2$ is the lower edge of the rectangular disc $D_{i,-}$.  \\

Noting that $X^{\prime \prime}_2$ and $X_1$ both represent $\cal{X}$ and have the same transverse knot type, we may find a transverse homeomorphism $g: \reals^3_- \rightarrow \reals^3_-$ which is the identity on $\reals^2$ with $g(X^{\prime \prime}_2) = X_1$.  Extend $g$ by the identity on $\reals^3_+$ to a transverse homeomorphism $G: \reals^3 \rightarrow \reals^3$.  Let $r_i = G(D_{i,-})$ and let $R_i = G(D_i) = r_i \cup D_{i, +}$.  The facts that $(1)$, $G$ is a transverse homeomorphism which is the identity on the half-space $\reals^3_+$ and that $(2)$ if $i \neq j$ then $D_{i,-} \cap D_{j,-} = D_i \cap D_j = \emptyset$, tell us that the $r_i$'s and $R_i$'s are pairwise disjoint embedded discs.  By Construction, $r_i \cap X_2 = \emptyset$,  whereas $R_i \cap X_2$ is a single point in the disc $D_{i,+}$.  Each $r_i$ joins $X^{\prime \prime}_2$ to $X_1$, meeting $X^{\prime \prime}_2$ in the arc $\beta_i$ and $X_1$ in the arc $\beta'_i = G(\beta_i)$.  Each $R_i$ joins $X'_2$ to $X_1$, meeting $X'_2$ in $\alpha_i$ and $X_1$ in $\beta_i$.  \\

Let $X'_3$ be the knot that is obtained from $X_1$ by replacing each $\beta_i \subset X_1$ by $\partial R_i - \beta'_i$.  Then $X'_3$ is constructed from $X_1$ by attaching $k$ pairwise disjoint long thin hooks to $X_1$.  There are two important aspects to this construction:

\begin{enumerate}
\item $X_2 \cup X'_3$ has the same transverse link type as $X_2 \cup X'_2$.  For, by construction, the transverse homeomorphism $G^{-1}: \reals^3 \rightarrow \reals^3$ sends $X_2 \cup X'_3$ to $X_2 \cup X'_2$.  
\item $X'_3$ is the connected sum of $X_1$ and $k$ copies of the unknot, the $i^{th}$ copy being $\partial R_i$.  
\end{enumerate}

Now recall from our initial construction that we had chosen $X'_2$ to be a transverse preferred longitude for $X_2$.  From $1.$ it follows that $X'_3$ is also a transverse preferred longitude for $X_2$.  (The linking number is unchanged).  There are two cases. \\

\np {\bf Case 1:} If $X'_3$ is transversely isotopic in the complement of the axis to a transverse closed braid, then let $X_3$ be that closed braid.  By construction, $X_3$ is the connected sum of $X_1$ and $k$ closed-braid copies of the unknot.  By the n-connect-sum Lemma, then, we conclude that $X_1 \equiv X_3$.  We already know that $X_3$ is a preferred longitude for $X_2$.  Changing our point of view, it follows that $X_2$ is a preferred longitude for $X_3$.  We conclude that $X_3 \equiv X_2$, by the TMT for preferred longitudes.  Therefore $X_1 \equiv X_2$.  \\

\np {\bf Case 2:} In general the transverse knot $X'_3$ will not be a closed braid because the arcs in $\partial r_i - \beta_i \cup \beta'_i$ will in general not be oriented to be increasing with respect to $\theta$.  In this case, we may apply Bennequin's Transverse Alexander Theorem to transversely change the subarcs that are wrongly oriented to a union of arcs that are correctly oriented and construct the transverse closed braid $X_3$.   In our situation, all of the subarcs which are not in braid position are in $\reals^3_-$, while $X_2$ is in $\reals^3_+$ , allowing us to assume that all such modifications take place away from $X_2$.  Therefore $X_2 \cup X_3$ has the same link type as $X_2 \cup X'_3$.  By the argument for Case $(1)$,  we have that $X_1 \equiv X_2$.  

\end{pf}

\end{document}